\title{Homological Connectivity of Random  \\ $k$-dimensional Complexes}
\documentstyle[amsmath,amsfonts,amssymb,amsthm,12pt]{article}
\font\sixbb=msbm6
\font\eightbb=msbm8
\font\twelvebb=msbm10 scaled 1095
\newfam\bbfam
\textfont\bbfam=\twelvebb \scriptfont\bbfam=\eightbb
                           \scriptscriptfont\bbfam=\sixbb
\def\bb{\fam\bbfam\twelvebb}

\newcommand{\Int}{{\bb Z}}

\newtheorem{theorem}{\bf Theorem}[section]
\newtheorem{claim}[theorem]{\bf Claim}

\newtheorem{proposition}[theorem]{\bf Proposition}

\newcommand{\enp}{\begin{flushright} $\Box$ \end{flushright}}
\newcommand{\beq}[0]{\begin{equation}}
\newcommand{\enq}[0]{\end{equation}}
\newcommand{\prob}{{\rm Pr}}

\newcommand{\supp}{{\rm Supp}}

\newcommand{\dn}{\Delta_{n-1}}

\newcommand{\cf}{{\cal F}}
\newcommand{\cg}{{\cal G}}

\newcommand{\spp}{{\rm supp}}

\begin{document}
\author{
{\sc R. Meshulam}\thanks
   {Institute for Advanced Study, Princeton, NJ 08540 and Department of Mathematics,
   Technion, Haifa 32000, Israel. Supported by a state of New Jersey grant  and by the
   Israel Science Foundation.
e-mail: meshulam@math.technion.ac.il} \and {\sc N. Wallach}\thanks
   {Department of Mathematics,
   Technion, Haifa 32000, Israel.  e-mail: tani@tx.technion.ac.il }}
%\date{}
%\begin{document}
%\insert\footins{\footnotesize\rule{0pt}{\footnotesep}
%\\ {\it Math Subject Classification.}  15A72, 05D05.
%\\ {\it Keywords and Phrases.} Exterior algebra, Skew-symmetric matrices, Matroid Parity \\ }
\maketitle
\pagestyle{plain}
\begin{abstract}
Let $\dn$ denote the $(n-1)$-dimensional simplex. Let $Y$ be a
random $k$-dimensional subcomplex of $\dn$ obtained by starting
with the full $(k-1)$-dimensional skeleton of $\dn$ and then adding
each $k$-simplex independently with probability $p$. Let
$H_{k-1}(Y;R)$ denote the $(k-1)$-dimensional reduced homology group of $Y$ with
coefficients in a finite abelian group $R$.
It is shown that for any fixed $R$ and $k \geq 1$ and for any function $\omega(n)$ that
tends to infinity
$$\lim_{n \rightarrow \infty} \prob ~[~ H_{k-1}(Y;R)=0 ~]= \left\{
\begin{array}{ll}
        0 & p=\frac{ k \log n-\omega(n)}{n} \\
        1 & p=\frac{ k \log n + \omega(n)}{n}
\end{array}
\right.~~ $$

\end{abstract}

\section{Introduction}
 Let $G(n,p)$ denote the probability space of graphs on the vertex
set $[n]=\{1,\ldots,n\}$ with independent edge probabilities $p$.
Let $\log$ denote the natural logarithm. A classical result of
Erd\H{o}s and R\'{e}nyi \cite{ER60} asserts that the threshold
probability for connectivity of $G \in G(n,p)$ coincides with the
threshold for the non-existence of isolated vertices in $G~$. In
particular, for any function $\omega(n)$ that tends to infinity
$$\lim_{n \rightarrow \infty} \prob ~[~G \in G(n,p)~:~G~ {\rm
~connected~}] =  \left\{
\begin{array}{ll}
        0 & p=\frac{\log n -\omega(n)}{n} \\
        1 & p=\frac{\log n+ \omega(n)}{n}
\end{array}
\right.~~.$$

A $2$-dimensional analogue of the Erd\H{o}s-R\'{e}nyi result was
considered in \cite{LM03}, where the threshold for homological
$1$-connectivity of random $2$-dimensional complexes was
determined (see below). In this paper we study the homological
$(k-1)$-connectivity of random $k$-dimensional complexes for a
general fixed $k$.

We recall some topological terminology (see e.g. \cite{Munkres}) .
Let $X$ be a finite simplicial complex on the vertex set $V$. let
$X^{(k)}=\{ \sigma \in X:\dim \sigma \leq k \}$ denote the
$k$-dimensional skeleton of $X$, and let $X(k)$ denote the set of
$k$-dimensional simplices in $X$, each taken with an arbitrary but
fixed orientation. Denote by $f_k(X)=|X(k)|$ the number of
$k$-dimensional simplices in $X$. Let $R$ be a fixed finite
abelian group of cardinality $r$. A simplicial $k$-cochain is an
$R$-valued skew-symmetric function on all ordered $k$-simplices of
$X$. For $k \geq 0$ let $C^k(X)$ denote the group of $k$-cochains
on $X$. The $i$-face of an ordered $(k+1)$-simplex
$\sigma=[v_0,\ldots,v_{k+1}]$ is the ordered $k$-simplex
$\sigma_i=[v_0,\ldots,\widehat{v_i},\ldots,v_{k+1}]$. The
coboundary operator $d_k:C^k(X) \rightarrow C^{k+1}(X)$ is given
by $$d_k \phi (\sigma)=\sum_{i=0}^{k+1} (-1)^i \phi
(\sigma_i)~~.$$ It is convenient to augment the cochain complex
$\{C^i(X)\}_{i=0}^{\infty}$ with the $(-1)$-degree term
$C^{-1}(X)=R$ with the coboundary map $d_{-1}:C^{-1}(X)
\rightarrow C^0(X)$ given by $d_{-1}a(v)=a$ for $a \in R~,~v \in
V$. Let $Z^k(X)= \ker (d_k)$ denote the space of $k$-cocycles and
let $B^k(X)={\rm Im}(d_{k-1})$ denote the space of
$k$-coboundaries. For $k \geq 0$ let $H^k(X;R)=Z^k(X)/B^k(X)~$
denote the $k$-th reduced cohomology group of $X$ with
coefficients in $R$. We abbreviate $H^k(X)=H^k(X;R)$.

Let $\dn$ denote the $(n-1)$-dimensional simplex on the vertex set
$V=[n]$.  Let $Y_k(n,p)$ denote the probability space of complexes
$\Delta_{n-1}^{(k-1)} \subset Y \subset \Delta_{n-1}^{(k)}$ with
probability measure
$$\Pr(Y)=p^{f_k(Y)}(1-p)^{\binom{n}{k+1}-f_k(Y)}~.$$
A $(k-1)$-simplex $\sigma \in \dn (k-1)$ is {\it isolated} in $Y$
if it is not contained in any of the $k$-simplices of $Y$. If
$\sigma$ is isolated then the indicator function of $\sigma$ is a
non-trivial $(k-1)$-cocycle of $Y$, hence $H^{k-1}(Y) \neq 0$. Our
main result is that the threshold probability for the vanishing of
$H^{k-1}(Y)$ coincides with the threshold for the non-existence of
isolated $(k-1)$-simplices in $Y$.
\begin{theorem}
\label{gen}  Let $k \geq 1$ and $R$ be fixed, and let $\omega(n)$ be any function which satisfies
$\omega(n) \rightarrow \infty$ then
\begin{equation}
\label{genp}\lim_{n \rightarrow \infty} \prob ~[~Y \in Y_k(n,p):
H^{k-1}(Y;R)=0 ~]= \left\{
\begin{array}{ll}
        0 & p=\frac{k\log n -\omega(n)}{n} \\
        1 & p=\frac{k\log n+ \omega(n)}{n}
\end{array}
\right.~~
\end{equation}
\end{theorem}
\noindent
{\bf Remarks:} \\
1. Theorem \ref{gen} remains true when $H^{k-1}(Y)$ is replaced by
the $(k-1)$-th reduced homology group $H_{k-1}(Y)=H_{k-1}(Y;R)$.
This follows from the universal coefficient theorem since
$H_{k-2}(Y)=0$ for $Y \in Y_k(n,p)$.
\\
2. The  $k=1$ case of Theorem \ref{gen} is the Erd\H{o}s-R\'{e}nyi
result. For $k=2$ and $R=\Int_2$ the theorem was proved in
\cite{LM03}. Our approach to the general case combines the method
of \cite{LM03} with some additional new ideas.
%In particular, the
%case $k=2$ of Proposition \ref{main} is a sharp version of
%Proposition 2.1 in \cite{LM03}. Proposition \ref{est1} which is
%the main technical tool in the present paper, is a $k$-dimensional
%strengthening of Proposition 2.3 in \cite{LM03}.
\ \\

The case $p=\frac{ k\log n-\omega(n)}{n}$ of Theorem \ref{gen} is
straightforward: Let $g(Y)$ denote the number of isolated
$(k-1)$-simplices of $Y$. Then $$E[g]=\binom{n}{k}(1-p)^{n-k}=
\Omega(\exp(\omega(n)))~~.$$ A standard second moment argument
then shows that
$$\prob[~H^{k-1}(Y) = 0~] \leq \prob [~g=0~] =o(1)~.$$

The case $p=\frac{ k\log n+\omega(n)}{n}$ is more involved. For a
$\phi \in C^{k-1}(\dn)$ denote by $[\phi]$ the image of $\phi$ in
$H^{k-1}(\dn^{(k-1)})$. Let $$b(\phi)=|\{\tau \in \dn(k)~
:~d_{k-1}\phi(\tau) \neq 0\}|~.$$ For any complex $Y \supset
\dn^{(k-1)}$ we identify $H^{k-1}(Y)$ with its image under the
natural injection $H^{k-1}(Y) \rightarrow H^{k-1}(\dn^{(k-1)})$.
It follows that for $\phi \in C^{k-1}(\dn)$
%\begin{equation}
%\label{probp}
$$
\prob[~ [\phi] \in
 H^{k-1}(Y)~]=(1-p)^{b(\phi)}~.
$$
%\end{equation}
For $\phi \in C^{k-1}(\dn)$ let $\spp(\phi)=\{ \sigma \in \dn(k-1): \phi(\sigma) \neq 0\}$. The
{\it weight} of such $\phi$ is defined by
$$w(\phi)= \min~\{~|\spp (\phi')|~:~\phi' \in C^{k-1}(\dn)~,~[\phi']=[\phi]~\}=
$$ $$ \min ~\{~|\supp (\phi+d_{k-2}\psi)|~:~ \psi \in
C^{k-2}(\dn)~\}.$$ A $k$-uniform hypergraph $\cf \subset \binom{[n]}{k}$
is {\it connected} if for any $\sigma,\tau \in \cf$
there exists a sequence $\sigma=\sigma_0,\ldots,\sigma_t=\tau \in
\cf$ such that $|\sigma_i \cap \sigma_{i-1}|=k-1$ for all $1 \leq
i \leq t$. Let
$$\cg_n=\{0 \neq \phi \in C^{k-1}(\dn)~:~\spp(\phi)~{\rm~is~connected~},~w(\phi)=|\spp(\phi)|\}~.$$
If $H^{k-1}(Y) \neq 0$ and $\phi \in C^{k-1}(\dn)$ is a cochain of
minimum support size such that $0 \neq [\phi] \in H^{k-1}(Y)$, then
$\phi \in \cg_n$. Therefore
$$\prob [~H^{k-1}(Y) \neq 0~] \leq \sum_{\phi \in \cg_n} \prob[~
[\phi] \in H^{k-1}(Y)~]= \sum_{\phi \in \cg_n} (1-p)^{b(\phi)}~.$$
Theorem \ref{gen} will thus follow from
\begin{theorem}
\label{esti} For $p=\frac{ k\log n+\omega(n)}{n}$
\begin{equation}
\label{estim} \sum_{\phi \in \cg_n} (1-p)^{b(\phi)} =o(1) ~~.
\end{equation}
\end{theorem}
The main ingredients in the proof of Theorem \ref{esti} are a
lower bound on $b(\phi)$ given in Section \ref{s:lbnd}, and an
estimate for the number of $\phi \in \cg_n$ with prescribed values
of $b(\phi)$ given in Section \ref{s:domination}. In Section
\ref{s:final} we combine these results to derive Theorem
\ref{esti}. The group $R$ and the dimension $k$ are fixed
throughout the paper. We use $c_i=c_i(r,k)$ to denote constants
depending on $r$ and $k$ alone.

\section{A lower bound on $b(\phi)$}
\label{s:lbnd} We bound $b(\phi)$ in terms of the weight $w(\phi)$.
\begin{proposition}
\label{main}
For $\phi \in C^{k-1}(\dn)$
\begin{equation}
\label{bw}
b(\phi) \geq \frac{n w(\phi)}{k+1}~~.
\end{equation}
\end{proposition}
\noindent {\bf Proof:} For an ordered simplex
$\tau=[v_0,\ldots,v_{\ell}]$ and a vertex $v \not\in \tau$, let
$v\tau=[v,v_0,\ldots,v_{\ell}]$. For $u \in V$ define $\phi_u \in
C^{k-2}(\dn)$ by
\begin{equation}
\label{cn} \phi_u(\tau) = \left\{
\begin{array}{ll}
        \phi(u \tau)  & u \not\in \tau \\
        0 & u \in \tau~~.
\end{array}
\right.
\end{equation}
Let $\sigma \in \dn(k-1)$ and $u \in V$. Then
$$
\phi(\sigma)-d_{k-2}\phi_u(\sigma)=
 \left\{
\begin{array}{ll}
    d_{k-1}\phi(u\sigma)  & u \not\in \sigma \\
        0 & u \in \sigma~~.
\end{array}
\right.
$$
It follows that $$ (k+1)|\spp(d_{k-1}\phi)|= |\{(\tau,u):u \in
\tau \in \spp(d_{k-1}\phi) \}|=
$$ $$
|\{(\sigma,u) \in \dn(k-1) \times V~:~\sigma \in
\spp(\phi-d_{k-2}\phi_u)\}| = $$ $$ \sum_{u \in V}
|\spp(\phi-d_{k-2}\phi_u)|  \geq n w(\phi)~~.$$ {\enp} \noindent
{\bf Remark:} The following example shows that equality can be
attained in (\ref{bw}). Let $n$ be divisible by $k+1$, and let
$[n]=\cup_{i=0}^k V_i$ be a partition of $[n]$ with
$|V_i|=\frac{n}{k+1}$. Consider the unique cochain $\phi \in
C^{k-1}(\dn)$ that satisfies
$$\phi([v_0,\ldots,v_{k-1}])=
\left\{
\begin{array}{ll}
        1  & v_i \in V_i {\rm~for~all~} 0 \leq i \leq k-1 \\
        0   &   |\{v_0,\ldots,v_{k-1}\} \cap
V_i| \neq 1 {\rm~for~some~ } 0 \leq i \leq k-1.
\end{array}
\right. $$ Then $b(\phi)=(\frac{n}{k+1})^{k+1}$, and it can be shown that
$w(\phi)=|\spp(\phi)|=(\frac{n}{k+1})^k.$

\section{The number of $\phi$ with prescribed $b(\phi)$}
\label{s:domination}

Let $$\cg_n(m)=\{\phi \in \cg_n: |\spp(\phi)|=m\}$$ and for $0
\leq \theta \leq 1$ let
$$\cg_n(m,\theta)=\{\phi \in \cg_n(m): b(\phi)=(1-\theta)mn\}~.$$
Write $g_n(m)=|\cg_n(m)|$ and $g_n(m,\theta)=|\cg_n(m,\theta)|$.
Proposition \ref{main} implies that $g_n(m,\theta)=0$ for $\theta
>\frac{k}{k+1}$.  Our main estimate is the following
\begin{proposition}
\label{est1} There exists a constant $c_1=c_1(r,k)$ such that for
any $n \geq 10k^2$, $m \geq
\frac{n}{2k}$, and $\theta \geq
\frac{1}{2k} $
\begin{equation}
\label{estt1} g_n(m,\theta) \leq \Bigl( c_1 \cdot
n^{(k-1)(1-\theta(1-\frac{1}{2k^2}))} \Bigr)^m~~.
\end{equation}
\end{proposition}
The proof of Proposition \ref{est1} depends on a certain partial
domination property of hypergraphs.
 Let $\cf \subset \binom{[n]}{k}$
be a $k$-uniform hypergraph of cardinality $|\cf|=m$. For $\sigma
\in \cf$ let
$$\beta_{\cf}(\sigma)= |\{ \tau \in \binom{[n]}{k+1}: \binom{\tau}{k}
\cap \cf =\{\sigma\} \}|~~$$ and let
$\beta(\cf)=\sum_{\sigma \in \cf} \beta_{\cf}(\sigma)$. Clearly
$\beta_{\cf}(\sigma) \leq n-k$ and $\beta (\cf) \leq m(n-k)$. For
$S \subset \cf$ let $$\Gamma(S)=\{\eta \in \cf: |\eta \cap \sigma|
= k-1 {\rm~for~some~} \sigma \in S\}~~.$$
\begin{claim}
\label{partdom} Let $0 < \epsilon \leq \frac{1}{2}$ and $n>2 \log
\frac{1}{\epsilon}+k$. Suppose that $$\beta(\cf) \leq
(1-\theta)m(n-k)$$ for some $0 < \theta \leq 1$. Then there exists
a subfamily $S \subset \cf$ such that
$$|\Gamma(S)| \geq (1-\epsilon) \theta m~~$$
and
$$|S| <
(20\log\frac{1}{\epsilon}) \cdot \frac{m}{n-k}
+2\log\frac{1}{\epsilon \theta}~~.$$
\end{claim}
\noindent {\bf proof:} Let $c(\epsilon)=2\log\frac{1}{\epsilon}$.
Choose a random subfamily $S \subset \cf$ by picking each $\sigma
\in \cf$ independently with probability
$\frac{c(\epsilon)}{n-k}~~.$ For any $\sigma \in \cf$ there exist
distinct $v_1, \ldots, v_{n-k-\beta_{\cf}(\sigma)} \in [n]
-\sigma$ and $\tau_1, \ldots, \tau_{n-k-\beta_{\cf}(\sigma)} \in
\binom{\sigma}{k-1}$ such that $\tau_i \cup \{v_i\} \in \cf$ for
all $i$. In particular
$$\prob[~\sigma \not\in \Gamma(S)~] \leq
\Bigl(1-\frac{c(\epsilon)}{n-k}\Bigr)^{n-k-\beta_{\cf}(\sigma)}~~,$$
hence
\begin{equation}
\label{leftout}
E[~|\cf -\Gamma(S)|~] \leq \sum_{\sigma \in \cf}\Bigl(1-\frac{c(\epsilon)}
{n-k}\Bigr)^{n-k-\beta_{\cf}(\sigma)}~~.
\end{equation}
Since $$\sum_{\sigma \in \cf}(n-k-\beta_{\cf}(\sigma)) =
m(n-k)-\beta(\cf)
\geq \theta m(n-k)$$
it follows by convexity from (\ref{leftout}) that
$$E[~|\cf -\Gamma(S)|~] \leq  (1-\theta)m + \theta m
\Bigl(1-\frac{c(\epsilon)}{n-k}\Bigr)^{n-k} \leq $$ $$ (1-\theta)m + \theta m
e^{-c(\epsilon)}=(1-\theta)m + \theta m \epsilon^2~~.$$
Therefore
$$E[~|\Gamma(S)|~] \geq (1-\epsilon^2) \theta m~~.
$$
Hence, since $|\Gamma(S)| \leq |\cf|=m$, it follows that
\begin{equation}
\label{bnd1} \prob [~|\Gamma(S)| \geq (1-\epsilon) \theta m ~] >
\epsilon(1-\epsilon)\theta~~.
\end{equation}
On the other hand
$$E[~|S|~]=\frac{c(\epsilon)m}{n-k}~~$$ and by the large deviation
inequality (see e.g. Theorem A.1.12 in
\cite{AS})
\begin{equation}
\label{bnd2} \prob[~|S|> \lambda \frac{c(\epsilon)m}{n-k}~] <
\Bigl(\frac{e}{\lambda}\Bigr)^{\lambda\frac{c(\epsilon)m}{n-k}}~~
\end{equation}
for all $\lambda \geq 1$.
Let
$$\lambda=10+\frac{n-k}{m}\Bigl(\frac{\log\frac{1}{\theta}}{\log\frac{1}{\epsilon}}+1\Bigr)$$
then
$$
\epsilon(1-\epsilon)\theta
>\Bigl(\frac{e}{\lambda}\Bigr)^{\lambda \frac{c(\epsilon)m}{n-k}}~.$$
Hence by (\ref{bnd1}) and (\ref{bnd2}) there exists an $S \subset
\cf$ such that $|\Gamma(S)| \geq (1-\epsilon)\theta m$ and
$$|S| \leq \lambda \frac{c(\epsilon)m}{n-k}=
(20\log\frac{1}{\epsilon}) \cdot \frac{m}{n-k}
+2\log\frac{1}{\epsilon \theta}~~.$$  {\enp} \noindent {\bf Proof
of Proposition \ref{est1}:} Define
$$\cf_n(m,\theta)=\{\cf \subset \binom{[n]}{k}~:~|\cf|=m~,~\beta(\cf) \leq
(1-\theta)mn \}$$ and let
$f_n(m,\theta)=\bigl|\cf_n(m,\theta)\bigr|$. If $\phi \in
\cg_n(m,\theta)$, then $\cf=\supp(\phi)\in \cf_n(m,\theta)$.
Indeed, if $\tau \in \binom{[n]}{k+1}$ satisfies $\binom{\tau}{k}
\cap \cf =\{\sigma\}$, then $d_{k-1}\phi(\tau)=\pm \phi(\sigma)
\neq 0$, hence $\beta(\cf) \leq b(\phi) =(1-\theta)mn$. Therefore
$$g_n(m,\theta) \leq (r-1)^m f_n(m,\theta)~~.$$ We next estimate
$f_n(m,\theta)$. Let $\cf \in \cf_n(m,\theta)$, then
$$\beta(\cf) \leq (1-\theta)mn =(1-\frac{\theta n-k}{n-k})m(n-k)~.$$ Applying Claim
\ref{partdom} with $\theta'= \frac{\theta n-k}{n-k}$ and
$\epsilon=\frac{1}{2k^2}$, it follows that there exists an $S
\subset \cf$ of cardinality $|S| \leq \frac{c_2 m}{n}$ with
$c_2=c_2(k)$, such that $|\Gamma(S)| \geq
(1-\frac{1}{2k^2})\theta' m$. The injectivity of the mapping
$$\cf \rightarrow (S,\Gamma(S),\cf-\Gamma(S))~~$$
implies that
$$ f_n(m,\theta) \leq \sum_{i=0}^{c_2m/n} \binom{\binom{n}{k}}{i} \cdot 2^{(\frac{c_2
m}{n})kn} \cdot \sum_{j=0}^{(1-\theta'(1-\frac{1}{2k^2}))m}
\binom{\binom{n}{k}} {j} \leq $$ $$c_3^m \binom{\binom{n}{k}}
{(1-\theta'(1-\frac{1}{2k^2}))m} \leq $$
$$
c_4^m
\Bigl(\frac{n^k}{m}\Bigr)^{(1-\theta'(1-\frac{1}{2k^2}))m}~~.$$
Therefore
$$g_n(m,\theta) \leq (r-1)^m f_n(m,\theta) \leq $$
$$ (r-1)^m c_4^m \Bigl(\frac{n^k}{m}\Bigr)^{(1-\theta'(1-\frac{1}{2k^2}))m} \leq
$$
$$
\Bigl(c_1 \cdot
n^{(k-1)(1-\theta(1-\frac{1}{2k^2}))}\Bigr)^m~~~.
$$
{\enp}

\section{Proof of Theorem \ref{esti}}
\label{s:final}

{\bf Proof of Theorem \ref{esti}:}  Let $\omega(n) \rightarrow
\infty$ and let $p=\frac{k\log n+\omega(n)}{n}$. We have to show
that
\begin{eqnarray}
\label{tm1} \sum_{m \geq 1} \sum_{\phi \in \cg_n(m)}
(1-p)^{b(\phi)}=o(1)~~.
\end{eqnarray}
We deal separately with two intervals of $m~$:\
\\
\\ (i) $~1 \leq m \leq \frac{n}{2k}~$.
If $\phi \in \cg_n(m)$ then $\spp(\phi) \subset \binom{[n]}{k}$
is a connected $k$-uniform hypergraph, hence there exists a subset $S \subset [n]$ of
cardinality $|S| \leq m+k-1$ such that $\spp(\phi) \subset \binom{S}{k}$.
Since $d_{k-1}\phi(u\sigma)=\phi(\sigma) \neq 0$ for any $\sigma \in \spp(\phi)$
and $u \not\in S$, it follows that
$b(\phi) \geq
m(n-m-k+1)$.
The trivial estimate
$$g_n(m) \leq (r-1)^m \binom{\binom{n}{k}}{m} \leq c_5^m \Bigl(\frac{n^k}{m}\Bigr)^m$$
implies that for $n \geq 6k$
$$g_n(m) (1-p)^{m(n-m-k+1)} \leq$$ $$
c_5^m  \frac{n^{km}}{m^m} \Bigl(1-\frac{k\log
n+w(n)}{n}\Bigr)^{m(n-m-k+1)} \leq $$
$$
c_5^m  \frac{n^{km}}{m^m}  n^{\frac{-k(n-m-k+1)m}{n}}
 e^{\frac{-w(n)(n-m-k+1)m}{n}} \leq $$
$$c_6^m  \Bigl(\frac{n^k}{m}  n^{\frac{-k(n-m)}{n}}\Bigr)^m
e^{-\frac{w(n)}{3}m} = $$ $$(c_6  \frac{n^{\frac{km}{n}}}{m}
e^{-\frac{w(n)}{3}})^m~~.$$ Since
$$
\frac{n^{\frac{km}{n}}}{m} \leq \left\{
\begin{array}{ll}
       n^{kn^{-1/3}}  & m \leq n^{2/3} \\
        n^{-1/6} & n^{2/3} \leq m \leq \frac{n}{2k}
\end{array}
\right.~~
$$
it follows that there exists a $c_7=c_7(r,k)$ such that for $m
\leq \frac{n}{2k}$ and $n\geq 6k$
$$g_n(m) (1-p)^{m(n-m-k+1)} \leq \Bigl( c_7 e^{-\frac{
w(n)}{3}} \Bigr)^m~~.$$  Therefore
$$\sum_{m=1}^{n/2k} \sum_{\phi \in \cg_n(m)}(1-p)^{b(\phi)} \leq
\sum_{m=1}^{n/2k} g_n(m)(1-p)^{m(n-m-k+1)}  \leq $$
\begin{equation}
\label{smallm} \sum_{m=1}^{n/2k} \Bigl( c_7 e^{-\frac{w(n)}{3}}
\Bigr)^m =O(e^{-\frac{w(n)}{3}})=o(1)~~.
\end{equation}
(ii)
$m \geq \frac{n}{2k}$.
%For each $m$, the number of $\theta$ such
%that $g_n(m,\theta) \neq 0$ is at most $n^{k+1}$.
Then
$$\sum_{m \geq n/2k} \sum_{\theta \leq 1/2k}
\sum_{\phi \in \cg_n(m,\theta)}(1-p)^{b(\phi)} =
$$ $$\sum_{m \geq n/2k} \sum_{\theta \leq 1/2k}
g_n(m,\theta)(1-p)^{(1-\theta)mn} \leq $$
$$
\sum_{m \geq n/2k} g_n(m)(1-p)^{(1-\frac{1}{2k})mn} \leq
$$
$$
\sum_{m \geq n/2k}
 \Bigl(\frac{c_5n^k}{m}\Bigr)^m n^{-(1-\frac{1}{2k})km} \leq $$
$$
\sum_{m \geq n/2k} (2k c_5 n^{k-1})^m n^{-(1- \frac{1}{2k})km} =
$$
\begin{equation}
\label{bigms} \sum_{m \geq n/2k} \Bigl(2k c_5 n^{-1/2}\Bigr)^m=
n^{-\Omega(n)}~~.
\end{equation}
Next note that by Proposition \ref{main}, $g_n(m,\theta)=0$ for
$\theta > \frac{k}{k+1}$. Hence, by Proposition \ref{est1}
$$
\sum_{m \geq n/2k} \sum_{\theta \geq 1/2k} \sum_{\phi \in
\cg_n(m,\theta)}(1-p)^{b(\phi)} =
$$ $$\sum_{m \geq n/2k} \sum_{\theta \geq 1/2k}
g_n(m,\theta)(1-p)^{(1-\theta)mn} \leq $$
$$
\sum_{m \geq n/2k} \sum_{\substack{ \theta \geq 1/2k     \\
g_n(m,\theta) \neq 0}}
 \Bigl( c_1 \cdot
n^{(k-1)(1-\theta(1-\frac{1}{2k^2}))} \Bigr)^m \cdot
n^{-(1-\theta)km} =$$ $$\sum_{m \geq n/2k}
\sum_{\substack{ \theta \geq 1/2k     \\
g_n(m,\theta) \neq 0}} \Bigl( c_1 \cdot
n^{\theta(1+\frac{k-1}{2k^2})-1} \Bigr)^m
 \leq
 n^{k+1}\sum_{m \geq n/2k}
 \Bigl( c_1 \cdot n^{\frac{k}{k+1}(1+\frac{k-1}{2k^2})-1} \Bigr)^m
=$$
\begin{equation}
\label{bigmb}
 n^{k+1}\sum_{m \geq n/2k}  \Bigl(c_1n^{-1/2k} \Bigr)^m=
n^{-\Omega(n)}~~.
\end{equation}
Finally (\ref{tm1}) follows from (\ref{smallm}), (\ref{bigms}) and
(\ref{bigmb}). {\enp}

\section{Concluding Remarks}
\label{s:cr}

We have shown that in the model $Y_k(n,p)$ of random $k$-complexes
on $n$ vertices,  the threshold for the vanishing of
$H^{k-1}(Y;R)$ occurs at $p=\frac{k \log n}{n}$, provided that
both $k$ and the finite coefficient group $R$ are fixed. One
natural concrete question is whether $p=\frac{k \log n}{n}$ is
also the threshold for the vanishing of $H^{k-1}(Y;\Int)~$.

More generally, in view of the detailed understanding of the
evolution of random graphs (see e.g. \cite{AS}), it would be
interesting to formulate and prove analogous statements concerning
the topology of random complexes. For example, what is the higher
dimensional counterpart of the remarkable double-jump phenomenon
that occurs in random graphs?

 \ \\ \\ {\bf ACKNOWLEDGMENT} \\ We would like
to thank Nati Linial for helpful discussions and comments.


\begin{thebibliography}{99}

\bibitem{AS}
N. Alon and J. Spencer, {\it The Probabilistic Method}, 2nd
Edition, Wiley-Intescience, 2000.

\bibitem{ER60}
P. Erd\H{o}s and A. R\'{e}nyi, On the evolution of random graphs,
{\it Publ. Math. Inst. Hungar. Acad. Sci.} {\bf 5}(1960) 17-61.

\bibitem{LM03}
N. Linial and R. Meshulam, Homological connectivity of random
2-complexes, {\it Combinatorica}, to appear.

\bibitem{Munkres}
J. Munkres, {\it Elements of Algebraic Topology}, Addison-Wesley,
1984.


\end{thebibliography}
\end{document}